\documentclass{article}

\usepackage[latin1]{inputenc}
\usepackage{newcent}
\usepackage{amssymb}
\usepackage{graphicx}

\title{Computation of Power Loss in
Likelihood Ratio Tests for Probability Densities Extended by Lehmann
Alternatives}
\author{Lucas Gallindo Martins Soares\\
\small\it
Departamento de Estatística e Informática\\
\small\it
Universidade Federal Rural de Pernambuco, Brasil\\
\normalsize\tt lucasgallindo@gmail.com}
\date{}
\begin{document}
\maketitle
\begin{abstract}
We compute the loss of power in likelihood ratio tests when we test
the original parameter of a probability density extended by the
first Lehmann alternative.
\end{abstract}
\section{Distributions Generated by Lehmann Alternatives}
In the context of parametric models for lifetime data,
\cite{GuptaQueSoAPorra} disseminated the study of distributions
generated by Lehmann alternatives, cumulative distributions that
take one of the following forms:
\begin{eqnarray}
G_{1}\left(x,\lambda\right)=\left[F(x)\right]^{\lambda}
&\mathrm{or}&G_{2}\left(x,\lambda\right)=1-\left[1-F(x)\right]^{\lambda}
\end{eqnarray}
where $F(x)$ is any cumulative distribution and $\lambda>0$. In the
present note, we are going to call both $G$ distributions
\emph{generated distributions} or \emph{extended distributions}. It
is easy to see that for integer values of $\lambda$, $G_{1}$ and
$G_{2}$ are, respectively, the distribution of the maximum and the
minimum of a sample of size $\lambda$, the support of the two
distribution is the same of $F$, and that the associated density
functions are
\begin{eqnarray}
g_{1}\left(x,\lambda\right)=\lambda
f(x)\left[F(x)\right]^{\lambda-1}
&\mathrm{and}&g_{2}\left(x,\lambda\right)=\lambda
f(x)\left[1-F(x)\right]^{\lambda-1}
\end{eqnarray}
where $f(x)$ is the density function associated with $F$. Suppose
that we generate a distribution $G(x|\lambda)$ based on the
distribution $F(x)$, and want to generate another distribution
$G'(x|\lambda,\lambda')$ repeating the process; It is easy to see
that the distribution $G'$ will be the same as $G$, for the new
parameter of the distribution, $\lambda\lambda'$ may be summarized
as a single one. This has the interesting side effect that the
standard uniparametric exponential distribution may be seen as a
distribution generated by the second Lehmann alternative from the
distribution $F(x)=1-e^{-x}$.

To compute the moments of distribution generated by Lehmann
alternatives, we use the change of variables $u=F(x)$ in the
expression
\begin{equation}
E\left[X^{k}|\lambda\right]=\int_{-\infty}^{\infty}x^{k}\lambda
f(x)\left[F(x)\right]^{\lambda-1}dx
\end{equation}
yielding
\begin{equation}
E\left[X^{k}|\lambda\right]=\int_{0}^{1}\lambda
Q^{k}(u)u^{\lambda-1}du=E_{\mathrm{Beta}(\lambda,1)}\left[Q(u)\right]
\end{equation}
where $Q(u)=F^{-1}(u)$ is the quantile function. This integral is
equivalent to the expectancy of $Q(u)$ with respect to a Beta
distribution with parameters $\alpha=\lambda,\beta=1$. The same
reasoning can be used to show that, for the second Lehmann
alternative,
$E\left[X^{k}|\lambda\right]=E_{\mathrm{Beta}(1,\lambda)}\left[Q(u)\right]$.

Using the log-likelihood functions
\begin{eqnarray}
G_{1}\left(x,\lambda\right)&=&n\ln\left(\lambda\right)+
\sum_{j=1}^{n}\ln
f\left(x_{j}\right)+\left(\lambda-1\right)\sum_{j=1}^{n}\ln
F\left(x_{j}\right)
\end{eqnarray}
and
\begin{eqnarray}
G_{2}\left(x,\lambda\right)&=&n\ln\left(\lambda\right)+
\sum_{j=1}^{n}\ln
f\left(x_{j}\right)+\left(\lambda-1\right)\sum_{j=1}^{n}\ln
\left[1-F\left(x_{j}\right)\right]
\end{eqnarray}
we see that the maximum likelihood estimators to the parameter
$\lambda$ have the forms
\begin{eqnarray}
\hat{\lambda}=-\frac{n}{\sum_{j=1}^{n}\ln F\left(x_{j}\right)}
&\mathrm{and}& \hat{\lambda}=-\frac{n}{\sum_{j=1}^{n}\ln
\left[1-F\left(x_{j}\right)\right]}
\end{eqnarray}

The existing literature about distributions generated by Lehmann
alternatives concerns mostly distributions defined on the interval
$(0,\infty)$ or in the real line, with the paper by
\cite{NadarajahKotz} being the more complete review of progresses
and the paper \cite{Nadarajah} being an interesting application of
the concepts developed outside the original proposal by
\cite{GuptaQueSoAPorra}, which was to analyze lifetime data. In the
present paper, we are concerned with some information theoretical
quantities of the first extension. These are not the only papers
dealing with the subject, but a complete list with comments would
be a paper on its own.
\section{Kullback-Leibler Divergence}
Given two probability density functions, the quantity defined as
\begin{eqnarray}
\label{Eq:KLDefinition}
D_{KL}\left(f|g\right)&=&\int_{\mathbb{R}}f(x)\ln\left(\frac{f(x)}{g(x)}\right)dx
\end{eqnarray}
is called Kullback-Leibler Divergence (abbreviated DKL) after the
authors of the classical paper \cite{KullbackLiebler}. Very often,
this quantity is used as a measure of distance between two
probability density functions, even though it is not a metric; This
divergence measure clearly is greater or equal than zero, with zero
occurring only and only if $f=g$, but it is not symmetric, so
$D_{KL}\left(f|g\right)\neq D_{KL}\left(g|f\right)$, and it does not
obey the triangle inequality also.

Rewriting equation (\ref{Eq:KLDefinition}), we get
\begin{eqnarray}
\int_{\mathbb{R}}f(x)\ln\left(\frac{f(x)}{g(x)}\right)dx&=&\int_{\mathbb{R}}f(x)\ln(f(x))-f(x)\ln(g(x))dx\\
&=&E_{f}\left[\ln(f(X))\right]-E_{f}\left[\ln(g(X))\right]
\end{eqnarray}
where $E_{f}\left[h(X)\right]$ is the expectation of the random
variable $h(X)$ with respect to the probability density $f$. Since
$D_{KL}\left(f|g\right)$ is greater than zero, we have that
\begin{equation}
E_{f}\left[\ln(f(X))\right]>E_{f}\left[\ln(g(X))\right]
\end{equation}
We will now show that maximizing the likelihood is equivalent to
minimize $D_{KL}\left(f|e\right)$, where $e$ is the empirical
distribution function. Calculating $D_{KL}\left(f|e\right)$ we
arrive at
\begin{eqnarray}
D_{KL}\left(f|e\right)&=&E_{f}\left[\ln(f(X))\right]-\frac{1}{n}\sum_{j=1}^{n}\ln\left(f(x_{j},\theta)\right)
\end{eqnarray}
where the rightmost term is the empirical log-likelihood multiplied
by a constant. So, maximizing the rightmost term we minimize the
whole divergence; Then the process of maximizing the likelihood is
equivalent to minimizing the divergence between the empirical
density and the parametric model. This result is very common in the
related literature, and is shown in full detail on sources like
\cite{EguchiCopas}, which gives an accessible but rather compact
deduction of properties of methods based on Likelihood Functions
using DKL. In the next (and last) section we draw freely from a
result shown in the \cite{EguchiCopas} paper that states that DKL
might be used to measure the loss of power in likelihood ratio tests
when the distribution under the alternative hypothesis is
mis-specified.

\section{Wrong Specification of Reference Distribution and Loss of Power in Likelihood Ratio Tests}
Suppose we have data from a probability distribution
$H(x|\theta,\lambda)$, and want to test the hypothesis that
$(\theta=\theta_0,\lambda=\lambda_0)$. The usual log-likelihood
ratio is expressed as
\begin{equation}
\Lambda(\lambda_0,\theta_0)=\frac{\ell(\hat{\lambda},\hat{\theta})}{\ell(\lambda_0,\theta_0)}
\end{equation}
where the notation $\hat{\xi}$ is used for the unrestricted maximum
likelihood estimative of the parameter $\xi$. Suppose we are not
willing to (or not able to) compute
$\ell(\hat{\lambda},\hat{\theta})$ because the estimative of the
parameter $\lambda$ is troublesome and decide to approximate the
likelihood ratio statistic using $\ell(\lambda_1,\tilde{\theta})$
instead of the likelihood under the alternative hypothesis, where
$\tilde{\theta}$ is the maximum likelihood estimator of $\theta$
given that $\lambda=\lambda_1$. We have then the relation
\begin{equation}
\Lambda(x)\approx\frac{\ell(\lambda_1,\tilde{\theta})}{\ell(\lambda_0,\theta_0)}
\end{equation}
A result by \cite{EguchiCopas}, section 3, states that the test
statistic generated this way is less powerful than the usual one,
with the loss in the power equal to
\begin{equation}
\Delta_{\mathrm{Power}}=
D_{KL}\left(f(x|\hat{\lambda},\hat{\theta}),f(x|\lambda_1,\tilde{\theta})\right)
\end{equation}

In the present paper, we are concerned with the case where the data
follows a distribution extended with the first Lehmann alternative,
where the original distribution is such that $F=F(x|\theta)$ for a
parameter $\theta$. The null hypothesis will be of the form
\begin{equation}
H_0:\theta=\theta_0,\lambda=1
\end{equation}
against a alternative hypothesis
\begin{equation}
H_A:\theta\neq\theta_0,\lambda\neq 1
\end{equation}
If we erroneously consider that the data doesn't come from a extended
distribution $G(x|\lambda,\theta)$, but from a population that follows the
original $F(x|\theta)$ distribution, we can say that we are approximating the
log-likelihood under the alternative hypothesis like in the previous discussion.
In this case, the log-likelihood will be taken under the hypothesis
\begin{equation}
H_{A'}:\theta\neq\theta_0,\lambda=1
\end{equation}
which generates the following expression for the log-likelihood:
\begin{equation}
\Lambda(x)\approx\frac{\ell(1,\tilde{\theta})}{\ell(1,\theta_0)}
\end{equation}
Then we have that the test has less power than the one using the full $G$ distribution;
The difference on the power of the tests is given by
\begin{equation}
\label{Eq:DeltaPower}
\Delta_{\mathrm{Power}}=
D_{KL}\left(g(x|\hat{\lambda},\hat{\theta})|g(x|1,\tilde{\theta})\right)
\end{equation}
The main point in the above discussion is that for testing
hypotheses about the "original" parameter $\xi$, the tests using the
extended version of distributions are always more powerful, with a
considerable difference in the error type II rate.

\begin{figure}[!hpb]
\caption{Loss of Power as a Function of $\lambda$, for $\lambda>1$.}
\label{Fig:PowerLossForBigLambda}
\includegraphics[angle=0,width=\textwidth]{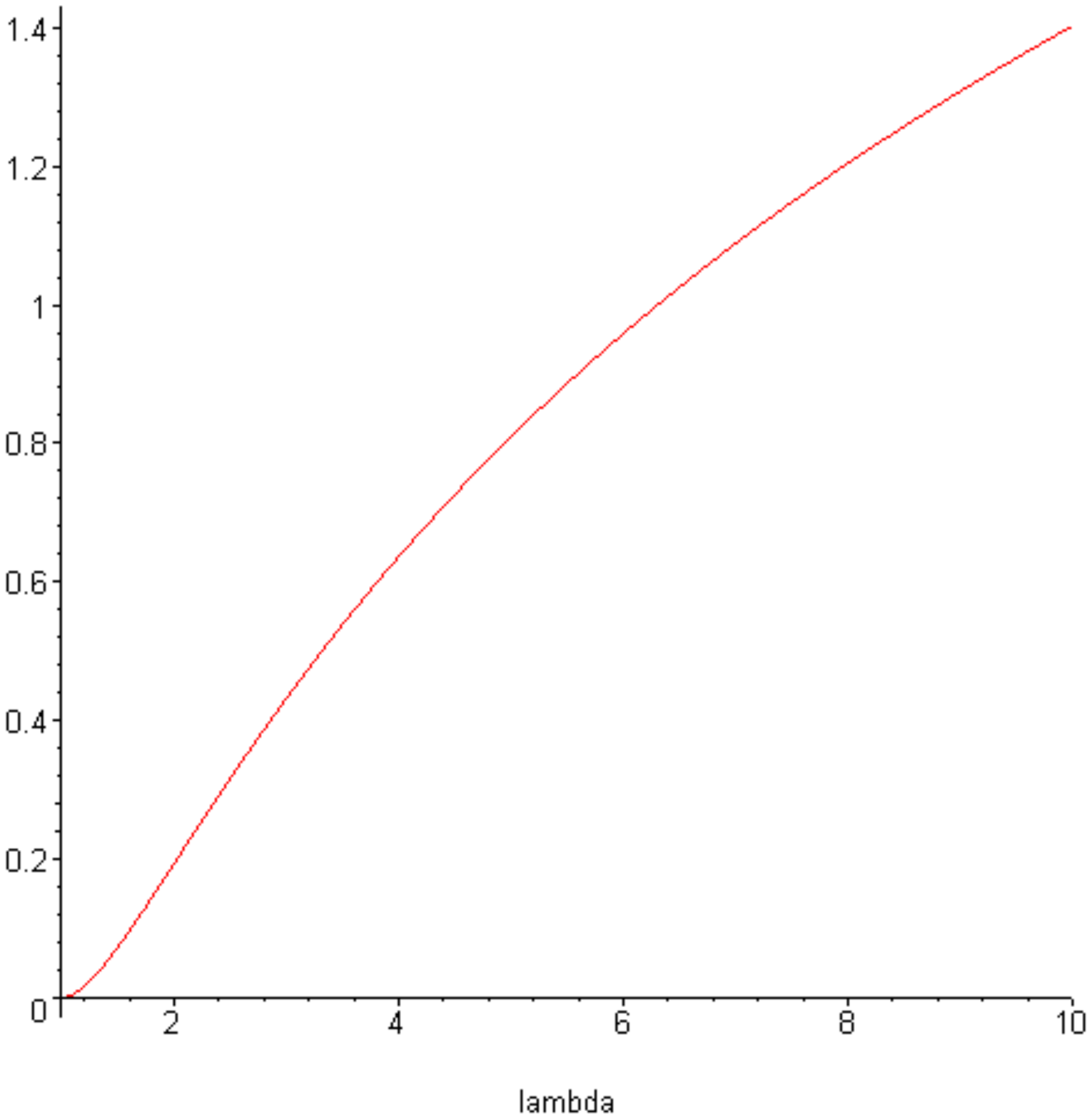}
\end{figure}

Expanding the equation (\ref{Eq:DeltaPower}) we have that
\begin{eqnarray}
\Delta_{\mathrm{P}}&=&D_{KL}\left(g(x|\hat{\lambda},\hat{\theta})|g(x|1,\tilde{\theta})\right)\\
&=&
\int_{\mathbb{R}}g(x|\hat{\lambda},\hat{\theta})\ln\left(\frac{g(x|\hat{\lambda},\hat{\theta})}{g(x|1,\tilde{\theta})}\right)dx\\
&=&\int_{\mathbb{R}}\lambda f(x|\hat{\lambda},\hat{\theta})F^{\lambda-1}(x|\hat{\lambda},\hat{\theta})\ln\left(\frac{\lambda f(x|\hat{\lambda},\hat{\theta})F^{\lambda-1}(x|\hat{\lambda},\hat{\theta})}{f(x|1,\tilde{\theta})}\right)dx\\
&=&\int_{\mathbb{R}}\lambda f(x|\hat{\lambda},\hat{\theta})F^{\lambda-1}(x|\hat{\lambda},\hat{\theta})\ln\left(\lambda F^{\lambda-1}(x|\hat{\lambda},\hat{\theta})\right)dx\\
&=&\ln\lambda+\int_{\mathbb{R}}\lambda(\lambda-1)f(x|\hat{\lambda},\hat{\theta})F^{\lambda-1}(x|\hat{\lambda},\hat{\theta})\ln\left(F(x|\hat{\lambda},\hat{\theta})\right)dx
\end{eqnarray}
Integrating by parts, we get
\begin{equation}
\label{Eq:DeltaPowerLambda}
\Delta_{\mathrm{Power}}=\ln\lambda+\frac{1-\lambda}{\lambda}
\end{equation}
The graphic of this function is the loss of power that we have
on our test when we the distribution of our data is one extended by the first Lehmann
alternative and we fail to notice that, and is depicted in Figure \ref{Fig:PowerLossForBigLambda}
for values of $\lambda$ bigger than one.

\begin{thebibliography}{}
\bibitem[Eguchi and Copas 1998]{EguchiCopas}\textsc{Eguchi, S. and Copas, J.} (2006).
Interpreting Kullback-Leibler divergence with the Neyman-Pearson
lemma. \emph{Journal of Multivariate Analysis}, vol. \textbf{97},
Issue 9, pages 2034-2040.
\bibitem[Gupta et alii 1998]{GuptaQueSoAPorra}\textsc{Gupta, R. C., Gupta, P. L. and Gupta, R. D.} (1998).
Modeling failure time data by Lehman alternatives.
\emph{Communication in Statistics: Theory and Methods}, vol.
\textbf{27}, pages 887-904.
\bibitem[Kullback and Leibler 1951]{KullbackLiebler}\textsc{Kullback, S. and Leibler, R. A.} (1951).
On information and sufficiency. \emph{The Annals of Mathematical
Statistics}, vol. \textbf{22}, Number 1, pages 79-86.
\bibitem[Nadarajah and Kotz 2006]{NadarajahKotz}\textsc{Nadarajah, S., Kotz, S.} (2006).
The Exponentiated Type Distributions. \emph{Acta Applicandae
Mathematicae}, vol. \textbf{92}, pages 97-111.
\bibitem[Nadarajah 2006]{Nadarajah}\textsc{Nadarajah, S.} 2006. The exponentiated Gumbel distribution
with climate application. \emph{Environmetrics}, vol. \textbf{17},
Number 1, pages 13-23.
\end{thebibliography}
\end{document}